\input amstex
\magnification=\magstep1 \baselineskip=13pt
\documentstyle{amsppt}
\vsize=8.7truein \CenteredTagsOnSplits \NoRunningHeads
\def\today{\ifcase\month\or
  January\or February\or March\or April\or May\or June\or
  July\or August\or September\or October\or November\or December\fi
  \space\number\day, \number\year}

\def\PP{{\bold P}}
\def\EE{{\bold E\thinspace }}

\def\per{\operatorname{per}}

\def\he{\operatorname{head}}
\def\ta{\operatorname{tail}}
\topmatter
\title Enumerating Contingency Tables via Random Permanents\endtitle
\author Alexander Barvinok \endauthor
\address Department of Mathematics, University of Michigan, Ann Arbor,
MI 48109-1043 \endaddress \email barvinok$\@$umich.edu
\endemail
\date March  2006 \enddate
\thanks This research was partially supported by NSF Grant DMS 0400617.
The author is grateful to Microsoft (Redmond) for hospitality
during his work on this paper.
\endthanks
\abstract Given $m$ positive integers $R=(r_i)$, $n$ positive
integers $C=(c_j)$ such that $\sum r_i = \sum c_j =N$, and $mn$
non-negative weights $W=(w_{ij})$, we consider the total weight
$T=T(R, C; W)$ of non-negative integer matrices (contingency tables)
$D=(d_{ij})$ with the row sums $r_i$, column sums $c_j$, and the
weight of $D$ equal to $\prod w_{ij}^{d_{ij}}$. We present a
randomized algorithm of a polynomial in $N$ complexity which
computes a number $T'=T'(R,C; W)$ such that 
$T' \leq T \leq \alpha(R, C) T'$ where 
$\alpha(R,C) = \min \left\{ \prod r_i! r_i^{-r_i}, \ \prod c_j! c_j^{-c_j} \right\} N^N/N!$.
In many cases, $\ln T'$ provides an asymptotically
accurate estimate of $\ln T$. The idea of the algorithm
is to express $T$ as the expectation of the permanent of
an $N \times N$ random matrix with exponentially distributed
entries and approximate the expectation by the integral $T'$ of an
efficiently computable log-concave function on ${\Bbb R}^{mn}$. Applications to counting 
integer flows in graphs are also discussed.
\endabstract
\endtopmatter
\document

\head 1. Introduction and main results \endhead

\subhead (1.1) Contingency tables \endsubhead Let us fix $m$
positive integers $r_1, \ldots, r_m$ and $n$ positive integers
$c_1, \ldots, c_n$ such that
$$r_1 + \ldots + r_m=c_1 + \ldots + c_n =N.$$
A non-negative $m \times n$ integer matrix $D=\left(d_{ij}\right)$
with the row sums $r_1, \ldots, r_m$ and the column sums $c_1,
\ldots, c_n$ is called a {\it contingency table} with the margins
$R=(r_1, \ldots, r_m)$ and $C=(c_1, \ldots, c_n)$. The problem of
efficient enumeration of contingency tables with prescribed
margins has attracted a lot of attention recently, see
\cite{DG95}, \cite{D+97}, \cite{CD03}, \cite{Mo02}, \cite{C+05}. The interest
in contingency tables is motivated by applications to statistics,
combinatorics and representation theory, cf. \cite{DG95} and
\cite{DG04}.

Let $W=\left(w_{ij}\right)$ be an $m \times n$ matrix of
non-negative weights $w_{ij}$. In this paper, we consider the
quantity
$$T(R, C; W)=\sum_D \prod_{ij} w_{ij}^{d_{ij}},$$
where the sum is taken over all contingency tables
$D=\left(d_{ij}\right)$ with the given margins $R=(r_1, \ldots,
r_m)$ and $C=(c_1, \ldots, c_n)$. Thus if $w_{ij}=1$ for all
$i,j$, the value of $T(R, C; W)$ is equal to the number of the
contingency tables with the given margins. If $w_{ij} \in
\{0,1\}$, the number $T(R, C; W)$ counts contingency tables $D$
for which we have $d_{ij}=0$ for all $i,j$ with $w_{ij}=0$ (here
we agree that $0^0=1$). In this case, $T(R,C; W)$ can be interpreted as the 
number of integer flows in a bipartite graph, see \cite{B+04} and \cite{C+05}. We note that counting 
integer flows in a general graph on $n$ vertices can be reduced to counting of integer flows 
in a bipartite graph on $n+n$ vertices and hence to counting weighted $n \times n$ 
contingency tables, see Section 1.5.  
 
Geometrically, one can view $T(R,C;W)$ as the generating function
over all integer points in the transportation polytope of $m
\times n$ non-negative matrices with the row sums $r_i$ and column
sums $c_j$, cf. \cite{BP99}.

We note that if $m=n$, $R=(1, \ldots, 1)$, and $C=(1,\ldots, 1)$
then
$$T(R,C; W)=\per W$$
is the {\it permanent} of the weight matrix $W$, that is,
$$\per W=\sum_{\pi} \prod_{i=1}^n w_{i \pi(i)},$$
where the sum is taken over all bijections $\pi: \{1, \ldots, n\}
\longrightarrow \{1, \ldots, n\}$, cf., for example, \cite{Mi78}.
A randomized polynomial time approximation algorithm to compute
the permanent of a given non-negative matrix was recently obtained
by M. Jerrum, A. Sinclair, and E. Vigoda \cite{J+04}.

We show that $T(R, C; W)$ can be represented as the expected
permanent of an $N \times N$ random matrix with exponentially
distributed entries.

Recall that a random variable $\gamma$ is {\it standard
exponential} if
$$\PP(\gamma \geq t)=\cases e^{-t} &\text{for\ } t > 0 \\ 1
&\text{for\ } t \leq 0. \endcases$$

Our starting point is the following result. \proclaim{(1.2)
Theorem} Given positive integers $r_1, \ldots, r_m$ and $c_1,
\ldots, c_n$ such that
$$r_1 + \ldots + r_m =c_1 + \ldots + c_n =N$$
and $mn$ real numbers $w_{ij}$, $i=1, \ldots, m$ and $j=1, \ldots,
n$, let us construct the $N \times N$ random matrix $A$ as
follows. The matrix $A=A(\gamma)$ is a function of the $m \times
n$ matrix $\gamma=\left(\gamma_{ij}\right)$ of independent
standard exponential random variables $\gamma_{ij}$. We represent
the set of rows of $A$ as a disjoint union of $m$ subsets $R_1,
\ldots, R_m$, where $|R_i|=r_i$ for $i=1, \ldots, m$ and the set
of columns of $A$ as a disjoint union of $n$ subsets $C_1, \ldots,
C_n$, where $|C_j|=c_j$. Thus $A$ is split into $mn$ blocks $R_i
\times C_j$. We sample $mn$ independent standard exponential
random variables $\gamma_{ij}$, $i=1, \ldots, m$ and $j=1, \ldots,
n$, and fill the entries of the block $R_i \times C_j$ of
$A=A(\gamma)$ by the copies $w_{ij} \gamma_{ij}$. Then the total
weight $T(R,C; W)$ of the $m \times n$ contingency tables with the
row sums $r_1, \ldots, r_m$ and column sums $c_1, \ldots, c_n$,
where the table $D=\left(D_{ij}\right)$ is counted with the weight
$$w(D)=\prod_{i,j} w_{ij}^{d_{ij}},$$
is equal to
$${\EE \per A \over r_1! \cdots r_m! c_1! \cdots c_n!}.$$
\endproclaim

We prove Theorem 1.2 in Section 2.

Although we can compute individual permanents $\per A$ via the
algorithm of \cite{J+04}, evaluating the expectation is still a
difficult problem. However, the expectation of an {\it approximate
permanent} of $A$ can be computed efficiently.

\subhead (1.3) An approximation algorithm to compute $T(R,C; W)$
\endsubhead
We rely heavily on the theory of {\it matrix scaling} and its
applications to approximating the permanent, in particular as
described in \cite{Lo71}, \cite{Si64}, \cite{KK96}, \cite{NR99}, \cite{L+00}, and
\cite{GS02}, as well as on the Markov chain based algorithms for
integrating log-concave densities \cite{AK91}, \cite{F+94}, 
\cite{FK99}, and \cite{Ve05}.

We assume here that the weights $w_{ij}$ are strictly positive,
which is not really restrictive since the zero weights can be
replaced by sufficiently small positive weights.

Let $A=\left(a_{ij}\right)$ be an $N \times N$ positive matrix.
Then there exist positive numbers $\xi_1, \ldots, \xi_N; \eta_1,
\ldots, \eta_N$ and a positive doubly stochastic (all row and
column sums are equal to 1) $N \times N$ matrix $B=B(A)$,
$B=\left(b_{ij}\right)$, such that
$$a_{ij}= b_{ij}\xi_i \eta_j \quad \text{for} \quad i,j=1, \ldots,
N.$$ Moreover, the matrix $B=B(A)$ is unique while the numbers
$\xi_1, \ldots, \xi_N$ and $\eta_1, \ldots, \eta_N$ are unique up
to a scaling $\xi_i \longmapsto \xi_i \tau$, $\eta_j \longmapsto
\eta_j \tau^{-1}$. This allows us to define the function
$$\sigma(A)=\prod_{i=1}^N \xi_i \eta_i.$$
We use the two crucial facts about $\sigma$:
\bigskip
$\bullet$ There is an algorithm, which, given a positive $N \times
N$ matrix $A$ and a number $\epsilon>0$ computes $\sigma(A)$
within relative error $\epsilon$ in time polynomial in $\ln
\epsilon^{-1}$ and $N$ \cite{L+00}
\medskip
and
\medskip
$\bullet$ The function $\sigma$ is {\it log-concave}, that is,
$$\ln \sigma\left(\alpha_1 A_1 + \alpha_2 A_2\right) \geq
\alpha_1 \ln \sigma\left(A_1\right) + \alpha_2 \ln \sigma
\left(A_2\right)$$ for any positive matrices $A_1$ and $A_2$ and
any non-negative $\alpha_1$ and $\alpha_2$ such that $\alpha_1
+\alpha_2=1$.
\bigskip
Our algorithm is based on replacing $\per A$ in Theorem 1.2 by the
(scaled) function $\sigma(A)$. Namely, we define
$$T'(R, C; W) ={N! \over N^N} {\EE \sigma(A) \over r_1! \cdots
r_m! c_1! \cdots c_n!}.$$ Since both $\sigma(A)$ and the
exponential density on ${\Bbb R}^{mn}$ are log-concave, and since
$\sigma(A)$ is efficiently computable for any positive $A$, we can
apply results of R. Kannan et al. \cite{AK91}, \cite{F+94}, and
\cite{FK99} and of L. Lov\'asz and S. Vempala \cite{Ve05}
 on efficient integration of log-concave functions to
show that there is a randomized fully polynomial time
approximation scheme to compute $T'(R, C; W)$.
\medskip
$\bullet$ We present a randomized algorithm, which, for any
$\epsilon
>0$ computes $T'(R, C; W)$ within relative error $\epsilon$ in
time polynomial in $\epsilon^{-1}$ and $N$ (in the unit cost
model).
\medskip

We discuss the details of our algorithm in Sections 3 and 4.
Namely, in Section 3 we present the necessary results regarding
$\sigma(A)$ while in Section 4 we discuss the integration problem.

Finally, we discuss how well the value of $T'(R,C; W)$
approximates $T(R,C;W)$.

\proclaim{(1.4) Theorem} For the number
$T'(R,C; W)$ computed by the algorithm of Section 1.3, we have
$$T'(R, C; W) \leq T(R, C; W) \leq \alpha(R, C) T'(R, C; W),$$
where 
$$\alpha(R, C)={N^N \over N!}\min\left\{ \prod_{i=1}^m {r_i! \over r_i^{r_i}}, \quad 
\prod_{j=1}^n {c_j! \over c_j^{c_j}} \right\}.$$
\endproclaim

We prove Theorem 1.4 in Section 5.

Let us consider the case of $m=n$ and
$$r_1= \ldots = r_n=c_1=\ldots=c_n=t.$$
Thus we are enumerating weighted {\it magic squares}, that is, $n
\times n$ square matrices with row and column sums equal to $t$.
Applying Stirling's formula in Theorem 1.4, we achieve the approximation factor of
$\alpha(R,C) \leq (nt)^{-1/2}(const \cdot t)^{n/2}$, that is, simply exponential in
the size $n$ of the matrix and polynomial in the line sum $t$.
Thus, for any fixed $t$, the algorithm of Section 1.3 can be
considered as an extension of the algorithm of N. Linial, A.
Samorodnitsky, and A. Wigderson \cite{L+00} for computing the
permanent of a positive matrix within a simply exponential factor.
On the other hand, if $n$ is fixed and $t$ grows, the number of
magic squares grows as a polynomial in $t$ of degree $(n-1)^2$,
see for example, \cite{St97}. Thus, in this case, the algorithm of
Section 1.3 allows us to capture the logarithmic order of $T(R,C;
W)$.

Apart from the case of $r_i=c_j=1$ (computation of the permanent),
most of the research thus far dealt with the case of $w_{ij}=1$,
that is, with the non-weighted enumeration of contingency tables.
M. Dyer, R. Kannan, and J. Mount \cite{D+97} showed that if the
margins are not too small, $r_i=\Omega\left(n^2m\right)$ and
$c_j=\Omega\left(m^2n\right)$, the Monte Carlo based approach
allows one to approximate the number of contingency tables within
a prescribed relative error $\epsilon>0$ in time polynomial in $m,
n$, and $\epsilon^{-1}$. In this case, the number of tables is
well approximated by the volume of the transportation polytope of
$m \times n$ non-negative matrices with the row sums $r_i$ and the
column sums $c_j$. Subsequently, B. Morris \cite{Mo02} improved
the bounds to $r_i=\Omega\left(n^{3/2} m \ln m \right)$ and
$c_j=\Omega\left(m^{3/2} n \ln n \right)$. The approximation we
get is much less precise, but applies to arbitrary weights
$W=\left(w_{ij}\right)$ and seems to be non-trivial even for
$w_{ij}=1$ and moderate values of $r_i, c_j$. For example, if
$m=n$ and $r_i=c_j=n$, we approximate the number of tables within
a factor of $(const \cdot n)^{(n-2)/2}$, while the exact number of
tables is at least $e^{O(n^2)}$. In other words, in many
non-trivial cases we get an asymptotically accurate estimate of
$\ln T(R, C; W)$.

Since every log-concave density can be arbitrarily closely
approximated by the push-forward (projection) of the Lebesgue
measure restricted to some higher dimensional convex body, the
algorithm of Section 1.3 can be viewed as a volume approximation
algorithm. In contrast to \cite{D+97} and \cite{Mo02}, the convex
body whose volume we approximate is not polyhedral.

\subhead (1.5) Counting integer flows in a graph \endsubhead
Let $G=(V,E)$ be a directed graph with the set $V$ of vertices and the set 
$E$ of edges. Hence every edge $e \in E$ is incident to the {\it head}  $\he(v) \in V$ of $e$ and 
the {\it tail} $\ta(e) \in V$. 
 We assume that $G$ is connected and that it does not contain loops or
multiple edges. Suppose further that each vertex $v$ has an integer number $a(v)$, called 
the {\it excess} of $v$, assigned to it, and that 
$$\sum_{v \in V} a(v)=0.$$
A set of non-negative integers $x(e): e \in E$ is called an {\it integer feasible flow} in $G$ if for every 
$v \in V$ the balance condition holds:
$$\sum_{e:\ head(e)=v} x(e) -\sum_{e: \ tail(e)=v} x(e)=a(v).$$    
If $G$ does not contain directed cycles $v_1 \rightarrow v_2 \rightarrow \ldots \rightarrow v_k
\rightarrow v_1$, the number of integer feasible flows is finite, possibly 0. 
The problem of efficient counting of integer feasible flows in a given graph has attracted some
attention recently, cf. \cite{B+04} and \cite{C+05}. A variation of the problem
involves introducing capacities of edges (upper bounds on the flows).

One can express the number of integer feasible flows in a graph with $|V|=n$ vertices as the number $T(R,C; W)$ of  weighted $n \times n$ contingency tables, where $w_{ij} \in \{0,1\}$ for 
all $i,j$. To this end, let us construct a bipartite graph with $n +n$ vertices as follows. For 
every vertex $v \in V$, we introduce the left copy $v_L$ and the right copy $v_R$. The 
directed edges  $u \rightarrow v$ of $G$ are represented by the edges $u_L \rightarrow v_R$ of the 
bipartite graph. We also introduce edges $v_L \rightarrow v_R$. Finally, let us choose 
a sufficiently large integer $z$, for example,
$$z=\sum_{v: \ a(v)>0} a(v)$$ and let us assign the excesses
$$a\left(v_L\right)=z-a(v) \quad \text{and} \quad a\left(v_R\right)=z.$$
With a feasible flow in the original graph $G$ we associate a feasible 
flow in the constructed bipartite graph by letting the flow on the edge 
$u_L \rightarrow v_R$ equal to the flow on the edge $u \rightarrow v$ and 
assigning the flow $v_L \rightarrow v_R$ so as to satisfy the 
balance conditions. This correspondence is a bijection between the integer 
feasible flows in $G$ and the bipartite graph. Hence the number of such flows 
is equal to the number of  weighted $n \times n$ contingency tables with the rows and 
columns indexed by the vertices $v \in V$, the row margins $z-a(v)$, the column margins 
$z$ and the matrix $W=(w_{ij})$ of weights defined by $w_{ij}=1$ for 
$(i,j) \in E$ and $w_{ij}=0$ for $(i,j) \notin E$.

\head 2. Proof of Theorem 1.2 \endhead

If $\gamma$ is a standard exponential random variable then for any
integer $d \geq 0$ we have
$$\EE \gamma^d =\int_0^{+\infty} \tau^d e^{-\tau} \ d \tau= d!.$$

Let us consider the random matrix $A=\left(a_{pq}\right)$ as
defined in Theorem 1.2. We identify both the set of rows of $A$
and the set of columns of $A$ with the set $\{1, \ldots, N\}$.

For every permutation $\pi: \{1, \ldots, N\} \longrightarrow \{1,
\ldots, N\}$, let
$$t_{\pi}=\prod_{k=1}^N a_{k \pi(k)} \tag2.1$$ be the corresponding term
of $\per A$. Thus
$$\EE \per A =\sum_{\pi} \EE t_{\pi},\tag2.2$$
where the sum is taken over all permutations $\pi$. With every
permutation $\pi$ we associate a contingency table $D=D(\pi)$,
called the {\it pattern} of $\pi$ as follows. We let
$D=\left(d_{ij}\right)$ where $d_{ij}$ is the number of indices $k
\in \{1, \ldots, N\}$ such that $k \in R_i$ and $\pi(k) \in C_j$,
so the $\bigl(k, \pi(k)\bigr)$th entry of $A$ lies in the block
$R_i \times C_j$ of $A$.

For the corresponding term $t_{\pi}$ of the permanent (2.1), we
have
$$\EE t_{\pi} = \prod_{ij}w_{ij}^{d_{ij}} d_{ij}!, \tag2.3$$
where $D=\left(d_{ij}\right)$ is the pattern of $\pi$.

Now, let us count how many permutations $\pi$ have a given pattern
$D=\left(d_{ij}\right)$. Let us represent each subset $R_i$ of
rows as a disjoint (ordered) union
$$R_i=\bigcup_{j=1}^n R_{ij} \quad \text{for} \quad i=1, \ldots, m$$ of
(possibly empty) subsets $R_{ij}$ with $|R_{ij}|=d_{ij}$ and each
subset $C_j$ of columns as a disjoint (ordered) union
$$C_j=\bigcup_{i=1}^m C_{ij} \quad \text{for} \quad j=1, \ldots, n$$
 of (possibly empty) subsets $C_{ij}$ with
$|C_{ij}|=d_{ij}$. This pair of partitions gives rise to exactly
$\prod_{ij} d_{ij}!$ permutation $\pi$ with the pattern $D$: we
choose $\pi$ in such a way that if $k \in R_{ij}$ then $\pi(k) \in
C_{ij}$ and we note that there are precisely $d_{ij}!$ bijections
$R_{ij} \longrightarrow C_{ij}$.

 On the other hand, the number
of partitions $R_i =\bigcup_{j} R_{ij}$ is
$${r_i! \over \prod_{j=1}^n d_{ij}!}$$
while the number of partitions $C_j= \bigcup_i C_{ij}$ is
$${c_j! \over \prod_{i=1}^m d_{ij}!}.$$
Therefore, the number of permutations with the given pattern
$D=(d_{ij})$ is
$${r_1! \cdots r_m! c_1! \cdots c_n! \over \prod_{ij} d_{ij}!}.$$
The proof now follows by (2.3) and (2.2). {\hfill \hfill \hfill}
\qed

\remark{Remark} Let us modify the definition of $A$ as follows:
instead of filling the $R_i \times C_j$ block by the copies of
$w_{ij} \gamma_{ij}$, we fill $R_i \times C_j$ by the copies of
just $w_{ij}$, so $A$ is constructed deterministically. It follows
from the proof above that the value of
$${\per A \over r_1! \cdots r_m! c_1! \cdots c_n!}$$
is equal to the total weight of the contingency tables with the
margins $r_1, \ldots, r_m$ and $c_1, \ldots, c_n$ provided the
weight of the table $D=\left(d_{ij}\right)$ is
$$\prod_{ij} {w_{ij}^{d_{ij}} \over d_{ij}!}$$
(the Fisher-Yates statistics, cf. \cite{DG95}).
\endremark
For another proof of Theorem 1.2 in a particular case of
$w_{ij}=1$, see \cite{Ba05}.

\head 3. Matrix scaling \endhead

Here we summarize the matrix scaling results that we need. All the
results in Theorem 3.1 below can be found in the literature

We reproduce the approach of L. Gurvits and A. Samorodnitsky
\cite{GS02} adapted to the case of the permanent (paper
\cite{GS02} treats a more general and more complicated setting of
mixed discriminants), which is, in turn, a modification of  D. London's \cite{Lo71} 
approach.

Also, we restrict ourselves to the case of
strictly positive matrices to avoid dealing with certain
combinatorial subtleties.

\proclaim{(3.1) Theorem} For every positive $N \times N$ matrix
$A=\left(a_{ij}\right)$ there exist unique positive $N$-vectors
$x=x(A)$, $y=y(A)$, and an $N \times N$ positive matrix $B=B(A)$
$$x=\left(\xi_1, \ldots, \xi_N\right), \quad y=\left(\eta_1, \ldots, \eta_N\right),
\quad \text{and} \quad B=\left(b_{ij}\right)$$ so that the
following holds
 \roster
\item We have
$$a_{ij} =b_{ij} \xi_i \eta_j \quad \text{for} \quad i,j=1, \ldots, N;$$
\item We have
$$\prod_{j=1}^N \eta_j=1;$$
\item Matrix $B$ is doubly stochastic, that is,
$$\sum_{i=1}^N b_{ij}=1 \quad \text{for} \quad j=1, \ldots, N \quad
\text{and} \quad \sum_{j=1}^N b_{ij}=1 \quad \text{for} \quad i=1,
\ldots, N.$$
\endroster
Let us define
$$\sigma(A)=\prod_{i=1}^N \xi_i, \quad \text{where}
\quad x(A)=(\xi_1, \ldots, \xi_N).$$ Then $\sigma$ is a
log-concave function on the set of positive matrices:
$$\ln \sigma\left( \alpha_1 A_1 + \alpha_2 A_2\right) \geq
\alpha_1 \ln \sigma \left( A_1 \right)+ \alpha_2 \ln \sigma
\left(A_2\right)$$ for any two positive $N \times N$ matrices
$A_1$ and $A_2$ and any two numbers $\alpha_1, \alpha_2 \geq 0$
such that $\alpha_1 + \alpha_2=1$.
\endproclaim
\demo{Proof} Let us consider the hyperplane
$$H=\left\{\left(\tau_1, \ldots, \tau_N\right): \quad \sum_{i=1}^N \tau_i=0
\right\}$$ in ${\Bbb R}^N$. With a positive matrix
$A=\left(a_{ij}\right)$, we associate the function $f_A: {\Bbb
R}^N \longrightarrow {\Bbb R}$,
$$f_A(t)=\sum_{i=1}^N \ln \left(\sum_{j=1}^N a_{ij} e^{\tau_j}
\right), \quad \text{where} \quad t=\left(\tau_1, \ldots,
\tau_N\right).$$ Then the restriction of $f_A(t)$ on $H$ is
strictly convex and, moreover, $f_A$ attains its unique minimum
$t^{\ast}=\left(\tau_1^{\ast}, \ldots, \tau_N^{\ast} \right)$,
$t^{\ast}=t^{\ast}(A)$, on $H$, see \cite{GS02}.

Since $f_A$ is smooth, $t^{\ast}$ is also a critical point and the
gradient of $f_A$ at the critical point is proportional to vector
$(1, \ldots, 1)$, from which we get
$$\sum_{i=1}^N \left({a_{ik} e^{\tau_k^{\ast}} \over \sum_{j=1}^N a_{ij}
e^{\tau_j^{\ast}}}\right)=\gamma \tag3.1.1$$ for some constant
$\gamma$ and $k=1, \ldots, N$.

Let
$$\xi_i = \sum_{j=1}^N a_{ij} e^{\tau_j^{\ast}} \quad \text{for}
\quad i=1, \ldots, N,$$ let
$$ \eta_j=e^{-\tau_j^{\ast}} \quad \text{for} \quad
j=1, \ldots, N,$$ and let us define an $N \times N$ matrix
$B=\left(b_{ij}\right)$ by
$$b_{ij}={a_{ij} \over \xi_i \eta_j} \quad \text{for} \quad
i,j=1, \ldots, N.$$ We note that
$$\prod_{j=1}^N \eta_j=1$$
since $t^{\ast}$ lies in the hyperplane $H$ with $\tau_1 + \ldots
+ \tau_N=0$.
 Then, by (3.1.1), we have
$$\sum_{i=1}^N b_{ij}= \gamma \quad \text{for} \quad j=1, \ldots, N
\tag3.1.2$$ On the other hand,
$$\sum_{j=1}^N b_{ij}=1 \quad \text{for}
\quad i=1, \ldots, N. \tag3.1.3$$ Since $B$ is a square matrix,
comparing (3.1.2) and (3.1.3), we infer that $\gamma=1$ and so we
established the existence of $x=(\xi_1, \ldots, \xi_N)$ and
$y=(\eta_1, \ldots, \eta_N)$ and $B$ satisfying (1)--(3).

To show uniqueness, we note that if $x=(\xi_1, \ldots, \xi_N)$,
$y=(\eta_1, \ldots, \eta_N)$, and $B$ satisfy (1)--(3), then we
must have
$$\xi_i =\sum_{j=1}^N a_{ij} \eta_j^{-1} \quad
\text{for} \quad i=1, \ldots, N$$ and hence, necessarily, the
point $t=(\tau_1, \ldots, \tau_N)$ defined by
$$\tau_j=-\ln \eta_j \quad \text{for} \quad j=1, \ldots, N$$
is a critical point of $f_A(t)$ on $H$. Since $f_A(t)$ is strictly
convex on $H$, there is a unique critical point
$t^{\ast}=t^{\ast}(A)$.

Thus function $\sigma(A)$ is well-defined. Moreover, we can write
$$\ln \sigma(A)= \sum_{i=1}^N \ln \xi_i =f_A\left(t^{\ast}\right)=
\min_{t \in H} f_A(t). \tag3.1.4$$ We observe that for any fixed
$t$, the function $g(A)=f_A(t)$ is concave on the set of positive
matrices $A=\left(a_{ij}\right)$.

Hence for any $t \in H$ and any $\alpha_1, \alpha_2 \geq 0$ such
that $\alpha_1 + \alpha_2 =1$, we have
$$f_{\alpha_1 A_1 + \alpha_2 A_2}(t) \geq
\alpha_1 f_{A_1}(t) + \alpha_2 f_{A_2}(t) \geq \alpha_1 \ln
\sigma\left(A_1\right) + \alpha_2 \ln \sigma \left(A_2 \right).$$
Taking the minimum over $t \in H$, we conclude that
$$\ln \sigma\left( \alpha_1 A_1 + \alpha_2 A_2 \right)
\geq \alpha_1\ln \sigma\left(A_1\right) + \alpha_2 \ln \sigma
\left(A_2 \right),$$ so $\sigma(A)$ is indeed log-concave. \qed
\enddemo

\remark{(3.2) Remark} Another useful property of $\sigma(A)$ which
easily follows from (3.1.4) is {\it monotonicity}: if
$A=\left(a_{ij}\right)$ and $A'=\left(a_{ij}'\right)$ are positive
matrices such that $a_{ij}' \leq a_{ij}$ for all $i$ and $j$ then
$\sigma(A') \leq \sigma(A)$. We also note that $\sigma(A)$ is
positive homogeneous of degree $N$: $\sigma(\lambda A)=\lambda^N
\sigma(A)$ for all positive $N \times N$ matrices $A$ and all
$\lambda>0$.
\endremark

\subhead (3.3) Computing $\sigma(A)$ \endsubhead N. Linial, A.
Samorodnitsky, and A. Wigderson present in \cite{L+00} a
deterministic polynomial time algorithm, which, given an $N \times
N$ positive matrix $A$ and a number $\epsilon>0$ computes the
value of $\sigma(A)$ within a factor of $(1+\epsilon)$ in time
polynomial in $\ln \epsilon^{-1}$ and $N$ (in the unit cost
model).

We are interested in computing $\sigma(A)$ where $A=A(\gamma)$ is
a random matrix of Theorem 1.2. Thus $A$ is positive with
probability 1. We observe that we can further save on computations
as follows.

Let us consider the $m \times n$ matrix $\left(w_{ij}
\gamma_{ij}\right)$, which is also positive with probability 1.
Applying the algorithm of \cite{L+00}, we can scale the matrix to
the row sums $r_i$ and the column sums $c_j$. Namely, we can
compute (approximately, in polynomial time) positive numbers
$\lambda_i$, $i=1, \ldots, m$, and $\mu_j$, $j=1, \ldots, n$, and
an $m \times n$ positive matrix $L=\left(l_{ij}\right)$ such that
$$w_{ij}\gamma_{ij}=l_{ij} \mu_i \lambda_j \quad
\text{for} \quad i=1, \ldots, m \quad \text{and} \quad j=1,\ldots,
n$$ and such that
$$\split &\sum_{j=1}^n l_{ij} = r_i \quad \text{for} \quad
i=1, \ldots, m \quad \text{and} \\ &\sum_{i=1}^m l_{ij} =c_j \quad
\text{for} \quad j=1, \ldots, n.
\endsplit$$

If we divide every row of $A$ from $R_i$ by $\mu_i r_i$ and divide
every column from $C_j$ by $\lambda_j c_j$, we get the $N \times
N$ matrix with the entries in the $R_i \times C_j$ block equal to
$l_{ij}/r_ic_j$. It is seen that the obtained matrix is doubly
stochastic. Therefore, we have
$$\sigma(A)=\left(\prod_{i=1}^m \left(\mu_i r_i \right)^{r_i}
\right) \left(\prod_{j=1}^n \left(\lambda_j c_j
\right)^{c_j}\right).$$ Hence the scaling of the $N \times N$
matrix $A$ reduces to the scaling of the $m \times n$ matrix
$\left(w_{ij} \gamma_{ij}\right)$.

\head 4. Integrating $\sigma(A)$ \endhead

Here we describe an algorithm for computing
$$T'(R,C; W)={N! \over N^N}{\EE \sigma(A) \over r_1! \cdots r_m! c_1! \cdots c_n!},$$
cf. Section 1.3.

\subhead (4.1) Notation \endsubhead We interpret the space ${\Bbb
R}^{mn}$ as the space of all $m \times n$ matrices
$\gamma=\left(\gamma_{ij}\right)$. Let ${\Bbb R}^{mn}_+$ denote
the positive orthant $\gamma_{ij} > 0$ of ${\Bbb R}^{mn}$ and let
$$\Delta=\Bigl\{\gamma: \quad \sum_{ij} \gamma_{ij}=1
\quad \text{and} \quad \gamma_{ij} > 0 \quad \text{for} \quad i=1,
\ldots, m \quad \text{and} \quad j=1, \ldots, n \Bigr\}$$ be the
standard (open) simplex in ${\Bbb R}^{mn}$.

For $0< \delta < 1/mn$ let us consider the $\delta$-{\it interior}
of $\Delta$:
$$\Delta_{\delta}=\Bigl\{\gamma: \quad \sum_{ij} \gamma_{ij}=1 \quad
\text{and} \quad \gamma_{ij} > \delta \quad \text{for} \quad i=1,
\ldots, m \quad \text{and} \quad j=1, \ldots, n \Bigr\}.$$
Geometrically, $\Delta_{\delta}$ is an open simplex lying strictly
inside $\Delta$.

For a $\tau \in {\Bbb R}$, let $\nu$ be the Lebesgue measure on
the affine hyperplane
$$\sum_{ij} \gamma_{ij}=\tau$$
induced by the Euclidean structure on ${\Bbb R}^{mn}$.

For a matrix $\gamma \in {\Bbb R}^{mn}_+$, let
$$P(\gamma)=\per A(\gamma) \quad \text{and let} \quad
S(\gamma)=\sigma\left(A(\gamma)\right),$$ where $A(\gamma)$ is the
matrix constructed in Theorem 1.2 and $\sigma$ is the function of
Theorem 3.1.

Thus we have
$$T(R, C; W)={1 \over r_1! \cdots r_m! c_1! \cdots c_n!}
\int_{{\Bbb R}^{mn}_+} P(\gamma) \exp\Bigl\{-\sum_{ij} \gamma_{ij}
\Bigr\} d \gamma$$ and
$$T'(R,C; W)={N! \over N^N} {1 \over r_1! \cdots r_m! c_1! \cdots c_n!}
\int_{{\Bbb R}^{mn}_+} S(\gamma) \exp\Bigl\{-\sum_{ij} \gamma_{ij}
\Bigr\} d \gamma,$$ where $d \gamma$ is the Lebesgue measure on
${\Bbb R}^{mn}$.

To apply the results of \cite{AK91}, \cite{F+94},  \cite{FK99} (see also  \cite{Ve05})
on efficient integration of log-concave functions, we modify the
problem to that of integration of $P(\gamma)$ and $S(\gamma)$
first on $\Delta$ and then on $\Delta_{\delta}$.

We use that both functions $P(\gamma)$ and $S(\gamma)$ are
positive homogeneous of degree $N$ and monotone on ${\Bbb
R}^{mn}_+$: if $\gamma=\left(\gamma_{ij}\right)$ and
$\gamma'=\left(\gamma_{ij}' \right)$ are positive matrices such
that
$$\gamma_{ij}' \leq \gamma_{ij} \quad \text{for all} \quad i,j,$$
then
$$P(\gamma') \leq P(\gamma) \quad \text{and} \quad S(\gamma') \leq
S(\gamma),$$ cf. Remark 3.2.

\proclaim{(4.2) Lemma} We have
$$\int_{{\Bbb R}^{mn}_+}
P(\gamma) \exp\Bigl\{-\sum_{ij} \gamma_{ij} \Bigr\} d \gamma=
{(N+mn-1)! \over \sqrt{mn}} \int_{\Delta} P(\gamma) \ d
\nu(\gamma).$$
\endproclaim
\demo{Proof} We note that $${\Bbb R}^{mn}_+=\bigcup_{\tau > 0}
\tau \Delta.$$

Since $d \nu \ d \tau=\sqrt{mn} \ d \gamma$, we get
$$\int_{{\Bbb R}^{mn}_+}
P(\gamma) \exp\Bigl\{-\sum_{ij} \gamma_{ij} \Bigr\} d \gamma= {1
\over \sqrt{mn}} \int_0^{+\infty} e^{-\tau}\left( \int_{\tau
\Delta} P(\gamma) \ d \nu(\gamma) \right) \ d \tau.$$ Since
$P(\gamma)$ is positive homogeneous of degree $N$, we conclude
that
$$\int_{\tau \Delta} P(\gamma) \ d \nu(\gamma)=\tau^{N+mn-1}
\int_{\Delta} P(\gamma) \ d \nu(\gamma),$$ from which the proof
follows. {\hfill \hfill \hfill} \qed
\enddemo

The same identity holds for the integrals of $S(\gamma)$.

Next, we approximate the integral over the simplex $\Delta$ by the
integral over the inner simplex $\Delta_{\delta}$.

\proclaim{(4.3) Lemma} Let $\delta \leq 1/mn$ be a non-negative
number. Then
$$\bigl(1-mn \delta\bigr)^{N+mn-1}
\int_{\Delta} P(\gamma) \ d \nu(\gamma) \leq
\int_{\Delta_{\delta}} P(\gamma) \ d \nu(\gamma) \leq
\int_{\Delta} P(\gamma) \ d \nu(\gamma).$$
\endproclaim
\demo{Proof} To prove the lower bound, we observe that the
transformation
$$\gamma_{ij} \longmapsto \gamma_{ij}-\delta \quad \text{for}
\quad i=1, \ldots, m \quad \text{and} \quad j=1, \ldots, n$$ maps
$\Delta_{\delta}$ inside $(1-\delta mn) \Delta$. Since $P$ is
monotone, we get
$$\int_{\Delta_{\delta}} P(\gamma) \ d \nu(\gamma) \geq
\int_{(1-\delta mn) \Delta} P(\gamma) \ d \nu(\gamma)= (1-\delta
mn)^{N+mn-1} \int_{\Delta} P(\gamma) \ d \nu(\gamma),$$ where we
used that $P$ is homogeneous of degree $N$.

The upper bound is obvious. {\hfill \hfill \hfill} \qed
\enddemo

The same inequalities hold for the integrals of $S(\gamma)$.

For an $0<\epsilon<1$, let us choose a positive
$$\delta \leq {-\ln(1-\epsilon) \over mn(N+mn-1)} \approx
{\epsilon \over mn(N+mn-1)} \quad \text{for small} \quad \epsilon
>0.$$
Then the integral
$$\int_{\Delta_{\delta}} P(\gamma) \ d \nu(\gamma)$$
approximates the integral
$$\int_{\Delta} P(\gamma) \ d \nu(\gamma)$$
within a factor of $(1-\epsilon)$ and the same holds for the
integrals of $S(\gamma)$.

Since $A(\gamma)$ depends linearly on $\gamma$, by the results of
Section 3, $S(\gamma)$ is a strictly positive log-concave function
on the set of positive matrices $\gamma$ and the value of
$S(\gamma)$ can be computed in polynomial time for any given
positive matrix $\gamma$.

Our goal consists of estimating $T(R, C; W)$ by
$$T'_{\delta}(R, C; W)={N! \over N^N} {1 \over r_1! \cdots r_m!
c_1! \cdots c_n!} {(N+mn-1)! \over \sqrt{mn}}
\int_{\Delta_{\delta}} S(\gamma) \ d \gamma.$$ To compute the
integral, we apply the algorithms of \cite{AK91}, \cite{F+94}, and
\cite{FK99}. The computational complexity of the algorithms is
polynomial in the dimension $mn-1$ of the integral {\it and} the
Lipschitz constant of $\ln S$ on $\Delta_{\delta}$. Hence it
remains to estimate the Lipschitz constant of $\ln S$.

\proclaim{(4.4) Lemma} Let $\delta < 1/mn$ be a positive number.
Let $\gamma=\left(\gamma_{ij}\right)$ and
$\gamma'=\left(\gamma_{ij}'\right)$ be two matrices such that
$$\gamma_{ij}, \gamma_{ij}' \geq \delta \quad \text{for all} \quad
i,j.$$ Then
$$|\ln S(\gamma)- \ln S(\gamma')| \leq {N \over \delta}
\max_{ij} | \gamma_{ij}-\gamma_{ij}'|.$$
\endproclaim
\demo{Proof} For $t=\left(\tau_1, \ldots, \tau_N\right)$, let
$$f_{\gamma}(t)=\sum_{p=1}^N \ln \left(\sum_{q=1}^N a_{pq}(\gamma) e^{\tau_q} \right),
\quad \text{where} \quad A(\gamma)=\left(a_{pq}(\gamma)\right)$$
is the matrix of Theorem 1.2. Letting
$$H=\left\{(\tau_1, \ldots, \tau_N): \quad \sum_{i=1}^N \tau_i =0
\right\},$$ by formula (3.1.4), we can write
$$\ln S(\gamma)= \min_{t \in H} f_{\gamma}(t). \tag4.4.1$$

Let
$$\alpha=\max_{ij} |\gamma_{ij}-\gamma_{ij}'|.$$
Then
$$\gamma_{ij} \leq \gamma_{ij}' + \alpha \leq
\gamma_{ij}'\left(1 + {\alpha \over \delta} \right) \quad
\text{for all} \quad i,j$$ and, similarly,
$$\gamma_{ij}' \leq \gamma_{ij} + \alpha \leq
\gamma_{ij}\left(1 + {\alpha \over \delta} \right) \quad \text{for
all} \quad i,j.$$

Since
$$a_{pq}(\gamma)=w_{ij} \gamma_{ij} \quad \text{provided} \quad
p \in R_i \quad \text{and} \quad q \in C_j,$$ we have
$$a_{pq}(\gamma) \leq a_{pq}(\gamma')\left(1 + {\alpha \over \delta}
\right)$$ and, similarly,
$$a_{pq}(\gamma') \leq a_{pq}(\gamma)\left(1 + {\alpha \over \delta}
\right).$$

Therefore, for all $t=\left(\tau_1, \ldots, \tau_N \right)$, we
have
$$f_{\gamma}(t) \leq f_{\gamma'}(t) + N \ln \left(1+{\alpha \over
\delta} \right) \leq f_{\gamma'}(t) + {\alpha N \over \delta}$$
and, similarly,
$$f_{\gamma'}(t) \leq f_{\gamma}(t) +{\alpha N \over \delta}.$$
Applying (4.4.1), we complete the proof. {\hfill \hfill \hfill}
\qed
\enddemo

Summarizing, we conclude that there is a randomized algorithm,
which, for any given $\epsilon >0$ computes the value of
$$T'(R, C; W)={N! \over N^N} {\EE \sigma(A) \over r_1! \cdots r_m! c_1! \cdots
c_n!},$$ where $A$ is the matrix of Theorem 1.2, within relative
error $\epsilon$ in time polynomial in $\epsilon^{-1}$ and $N$ (in
the unit cost model).

\head 5. Proof of Theorem 1.4 \endhead

Our proof is based on two estimates for the permanent of a
non-negative matrix.

\subhead (5.1) The van der Waerden bound \endsubhead Let
$B=\left(b_{ij}\right)$ be an $N \times N$ doubly stochastic
matrix, that is, a non-negative matrix with all row and column
sums equal to 1. Then
$$\per B \geq {N! \over N^N}.$$
This bound constituted van B.L. der Waerden's conjecture proved by
G.P. Egorychev \cite{Eg81} and D.I. Falikman \cite{Fa81}, see also
Chapter 12 of \cite{LW01}.

\subhead (5.2) A continuous extension of the Minc-Bregman bound
\endsubhead
Let $B=\left(b_{ij}\right)$ be an $N \times N$ non-negative
matrix. Let
$$s_i=\sum_{j=1}^N b_{ij} \quad \text{for} \quad i=1, \ldots, N$$
be the row sums of $B$.

If $b_{ij} \in \{0,1\}$, the bound
$$\per B \leq \prod_{i=1}^N (s_i!)^{1/s_i}$$
was conjectured by H. Minc and proved by L.M. Bregman \cite{Br73},
see also Chapter 2 of \cite{AS00}.

A. Samorodnitsky communicated to the author the following extension of the Minc-Bregman 
bound. Suppose that 
$$\sum_{j=1}^N b_{ij}=1 \quad \text{for all} \quad \text{for} \quad i=1, \ldots, N \tag5.2.1$$
and that 
$$b_{ij} \leq 1/t_i \quad \text{for all} \quad i=1, \ldots, N, \tag5.2.2$$
where $t_i, i=1, \ldots, N$, are positive integers. 
Then 
$$\per B \leq \prod_{i=1}^N {\left(t_i!\right)^{1/t_i}  \over   t_i}. \tag5.2.3$$
To deduce (5.2.3), we argue that the maximum of $\per B$ on the class of $N \times N$ non-negative matrices 
satisfying (5.2.1) and (5.2.2) is attained at a matrix with $b_{ij} \in \{0, 1/t_i\}$ for all $i,j$.
Indeed, let us choose a particular row index $i$. Then any non-negative matrix $B$ satisfying 
(5.2.1)--(5.2.2) can be written as a convex combination of two non-negative matrices $B'$ and $B''$ which satisfy 
 (5.2.1)--(5.2.2), agree with $B$ in all rows, except possibly the $i$th row, and, additionally,  satisfy 
 $b'_{ij}, b_{ij}'' \in \{0, 1/t_i \}$. Since the function $\per B$ is linear in every row, we conclude that 
 $\per B \leq \max\left\{ \per(B'), \per(B'') \right\}$. Proceeding as above for rows 
 $i=1, \ldots, N$, we may assume that $b_{ij} \in \{0, 1/t_i\}$ for all $i,j$, so (5.2.3) follows from the 
 Minc-Bregman bound.

A similar bound is obtained by G.W. Soules  \cite{So03}. If $B$ is a non-negative matrix satisfying 
(5.2.1)--(5.2.2) where $t_i$ do not have to be integer, then 
$$\per B \leq \prod_{i=1}^N {\Gamma^{1/t_i}(t_i+1) \over t_i}.$$

\demo{Proof of Theorem 1.4} For a given $m \times n$ positive
matrix $\gamma=(\gamma_{ij})$, let $A=A(\gamma)$ be the matrix
constructed in Theorem 1.2 and let $B=B(A)$ be the matrix
constructed in Theorem 3.1. We have
$$\per A= \sigma(A) \per B,$$
and we estimate $\per B$.

By the estimate of Section 5.1, we have
$$\per A \geq {N! \over N^N} \sigma(A),$$
from which
$$T(R,C;W) \geq T'(R,C; W).$$
On the other hand, as is discussed in Section 3.3, we can
construct $B=B(\gamma)$ as follows: first, we construct a positive
$m \times n$ matrix $L=\left(l_{ij}\right)$ such that
$$w_{ij} \gamma_{ij} =l_{ij} \mu_i \lambda_j \quad \text{for all}
\quad i,j$$ and some positive numbers $\mu_1, \ldots, \mu_m$ and
$\lambda_1, \ldots, n$ and such that
$$\sum_{j=1}^n l_{ij}=r_i \quad \text{for} \quad i=1, \ldots, m$$
and
$$\sum_{i=1}^m l_{ij}=c_j \quad \text{for} \quad j=1, \ldots, n$$
and then fill the $R_i \times C_j$ block of $B$ by $l_{ij}/r_i
c_j$.

It follows then that every entry of $B$ in the block $R_i$ of rows  does not exceed $1/r_i$.
Applying the bound of Section 5.2, we get
$$\per B \leq \prod_{i=1}^m {r_i! \over  r_i^{r_i}}  .$$
Similarly, every entry of $B$ in the block $C_j$ of columns does not exceed $1/c_j$, so we 
get 
$$\per B \leq \prod_{i=1}^n {c_j! \over  c_j^{c_j}}.$$

Therefore,
$$\per A \leq  \sigma(A) \min \left\{\prod_{i=1}^m { r_i!  \over r_i^{r_i}}, \quad  \prod_{i=1}^n 
{c_j! \over c_j^{c_j}}
\right\}. $$
and the proof follows. {\hfill \hfill
\hfill} \qed
\enddemo

\head Acknowledgments \endhead

The author is grateful to Leonid Gurvits for useful conversations and references and 
to Alex Samorodnitsky whose communication \cite{Sa06} resulted in sharper bounds in
Theorem  1.4.

\enddocument
\end